\documentclass[a4paper,twoside,12pt]{amsart}
\usepackage{amssymb,amscd}
\usepackage{hyperref}
\usepackage{color}
\theoremstyle{plain}
\newtheorem{theorem}{Theorem}
\newtheorem{example}[theorem]{Example}
\newtheorem{definition}[theorem]{Definition}

\newtheorem{cor}[theorem]{Corollary}
\newtheorem{rmk}[theorem]{Remark}
\newtheorem{lemma}[theorem]{Lemma}
\newtheorem{ass}[theorem]{Assumption}
\newtheorem{prop}[theorem]{Proposition}

\newtheorem*{rmq}{Remark}
\newtheorem*{lemma*}{Lemma}
\newtheorem*{claim}{Claim}
\newtheorem*{observation}{Observation}

\newtheorem*{theorem*}{Theorem}
\newcommand{\bP}{{\mathbb{P}}}
\newcommand{\C}{{\mathbb{C}}}

\newcommand{\Q}{{\mathbb{Q}}}

\newcommand{\Z}{{\mathbb{Z}}}

\newcommand{\OO}{{\mathcal O}}

\newcommand{\into}{\hookrightarrow}
 \newcommand{\rarrow}[3]{\smash{\mathop{\hbox to#3{\rightarrowfill}}\limits
^{\scriptstyle#1}_{\scriptstyle#2}}}
 \newcommand{\larrow}[3]{\smash{\mathop{\hbox to#3{\leftarrowfill}}\limits
^{\scriptstyle#1}_{\scriptstyle#2}}}
 \newcommand{\darrow}[3]{\llap{$\scriptstyle #1$}
\left\downarrow\vbox to#3{}\right.\rlap{$\scriptstyle #2$}}
 \newcommand{\uarrow}[3]{\llap{$\scriptstyle #1$}
\left\uparrow\vbox to#3{}\right.\rlap{$\scriptstyle #2$}}
  \newcommand{\Vline}[3]{\llap{$\scriptstyle #1$}
\left\Vert\vbox to#3{}\right.\rlap{$\scriptstyle #2$}}
\newcommand{\mapdown}[1]{\Big\downarrow
         \rlap{$\vcenter{\hbox{$\scriptstyle#1$}}$}}
\newcommand{\mapup}[1]{\Big\uparrow
         \rlap{$\vcenter{\hbox{$\scriptstyle#1$}}$}}
\newcommand{\mapleft}[1]{\mathop{\vbox{\ialign{
                                ##\crcr
    ${\scriptstyle\hfil\;\;#1\;\;\hfil}$\crcr
    \noalign{\kern2pt\nointerlineskip}
    \leftarrowfill\crcr}}\;}}
\newcommand{\mapright}[1]{{\vbox{\ialign{
                                ##\crcr
    ${\scriptstyle\hfil\;\;#1\;\;\hfil}$\crcr
 \noalign{\kern2pt\nointerlineskip}
    \rightarrowfill\crcr}}\;}}

\newcommand{\End}{\mathop{\rm End}\nolimits}

\newcommand{\alb}{\mathop{\rm Alb}\nolimits}
\newcommand{\alg}{{\rm alg}}

\newcommand{\blow}{\mathop{\rm Bl}\nolimits}

\newcommand{\chowalg}{\chow^{\rm alg}}

\newcommand{\chowchom}{\chow_{\rm hom}}
\newcommand{\chowhom}{\chow^{\rm hom}}
\newcommand\chowlin{\chow^{\rm lin}}
\newcommand{\chow}{\mathop{\rm CH}\nolimits}
\newcommand{\cl}{\mathop{\rm cl}\nolimits}

\newcommand{\comp}{{\scriptstyle \circ}}
\newcommand{\corr}{\mathop{\rm Corr}\nolimits}

\newcommand{\id}{\mathop{\rm id}\nolimits}

\newcommand{\pic}{\mathop{\rm Pic}\nolimits}
\newcommand{\rset}{\ \}}  % for }
\newcommand{\lset}{\  \{} % for { 
\newcommand{\set}[1]{\lset #1 \rset} 
\newcommand{\sett}[2]{\lset #1 \mid  #2 \rset} 

\newcommand{\supp}{\mathop{\rm supp}\nolimits}
\newcommand{\tr}{{}^{\mathsf T}\kern-0.9pt} %transpose: use for "b" 

\title[Bloch-type conjectures for threefolds]{Bloch-type conjectures and an example of a threefold of general type}

\date{January 2009}
\author{Chris PETERS}
\address{Department of Mathematics,  University of Grenoble I, 
UMR 5582 CNRS-UJF, \\
38402-Saint-Martin d'H{\`e}res, France
}
\email{chris.peters@ujf-grenoble.fr}

\begin{document}
\begin{abstract}  The hypothetical existence of a good theory of mixed motives
predicts many deep phenomena related to algebraic cycles. One of these, a
generalization of Bloch's conjecture says that ``small Hodge diamonds'' go
with ``small Chow groups''.  
Voisin's method \cite{V} (which produces examples with small Chow groups) is  analyzed carefully to widen its applicability.  A threefold of general
type without $1$- and $2$-forms is exhibited  for which this  extension   yields  Bloch's generalized
conjecture.
\end{abstract}
\maketitle

\vskip 5 pt
\centerline{\tt MSC Classification 14C15, 14C30}
\vskip 10 pt

\section{Introduction}
Very little is known about the Chow groups of algebraic varieties. This is 
even true for $0$-cycles on surfaces. Mumford \cite{Mum} has shown that if Albanese equivalence and rational equivalence coincide on $0$-cycles of degree $0$, then there are no holomorphic
$2$-forms. In this case the group of $0$-cycles of degree $0$ modulo rational equivalence
is \emph{representable}, i.e. isomorphic to an abelian variety (the Albanese
variety in this case). Bloch conjectured the converse and in \cite{B-K-L}
 this is shown for surfaces which are not of general type. Later it was
shown for several classes of surfaces of general type without  holomorphic
$2$-forms (see \cite{I-M}, \cite{B}, \cite{V}). For higher dimensions an
analogue of Mumford's result can be found in \cite{B-S} and \cite{E-L}.

An analogue of (a weaker form of) Bloch's conjecture can be formulated in terms
of  the \emph{level} of a non-zero Hodge structure $H=\bigoplus H^{p,q}$, i.~e.
the largest difference $|p-q|$ for which $H^{p,q}\not=0$. In particular, for a
projective manifold $X$ the cohomology group $H^k(X,\C)$ has level $< (k-2p)$
if and only if $H^{k,0}(X)=\dots= H^{k-p,k}(X)=0$. For instance $H^1$ has level
$<1$ if it vanishes and $H^2$ has level $<1$ if $H^2=H^{1,1}$. For surfaces
the two combined mean  exactly that $q=p_g=0$. Bloch's conjecture in this case
states that for $0$-cycles rational and homological equivalence coincide.  Laterveer   \cite[Cor.~1.10]{La}  and C.  Schoen \cite{Scho} have found a  generalization of Mumford's theorem:  if for all cycles of dimension $\le s$ rational and homological equivalence coincide, then every cohomology group  $H^k(X)$ has level  $< k-2s$.
Bloch's generalized conjecture is the converse. See also \cite{Ja},
where such conjectures are deduced from the Beilinson's conjectures on the
existence of a nice filtration on the Chow groups.  

An obvious class to test these types of conjectures are hypersurfaces of degree
$d$ in projective space $\bP^{n+1}$. By Lefschetz hyperplane theorem only the 
middle cohomology group can have level $>0$. In fact, it has level $<n-2s$ 
if $(s+1)d <n+2 $ and so one expects rational and homological equivalence
to coincide for cycles of dimension $\le s$. 

For $s=0$ this is true (Roitman's theorem \cite{R1}), but for $s>0$ the
Bloch conjecture  has only been shown for degrees much below the optimum which
is the largest  integer  $<(n+2)/(s+1)$.  See \cite{E-L-V} and
\S~\ref{sec:linchow} below.  Note that   $d$  falls below the optimum if one
lets $n$ grow. The same is true for the bound from \cite{E-L-V}.  Note also that
such hypersurfaces are Fano.

To get examples of projective varieties which on the one hand have small Hodge
diamonds, but on the other hand have ample canonical bundle, one has to look
for complete intersections in weighted projective spaces. The disadvantage is
that one has to allow ``quasi-smooth'' complete intersections which have certain
mild singularities. But, provided one works with rational
coefficients the Chow groups behave as in the smooth case. Also, the
cohomology groups carry pure Hodge structures as in the smooth case. See \S
\ref{subsec:wci} and the references therein for details. These examples thus
provide more interesting testing ground for Bloch-type conjectures since they
are far from Fano.

 Of crucial importance here is that the analogs of the bounds
from \cite{E-L-V} are valid for weighted projective spaces
$P=\bP(a_0,\dots,a_{n+1})$ with $n$ replaced by the sum of the weights $\sum
a_i:=N+1$. See \S~\ref{sec:linchow}. In other words, a complete intersection in
$P$  whose multidegree   is ``small'' with respect to $N$ has ``small''
Chow groups\footnote{The bounds one  gets in this way are again not the optimal ones in general.}. 
 
In  \cite{V}  Voisin considers  a variant of the Bloch conjecture  for certain smooth hypersurfaces in $\bP^3$ and $\bP^4$ equipped with a finite group of automorphisms;  instead of looking to the entire cohomology and the entire Chow group, one restricts to  one or more character subspaces at a time. In loc. cit.  this modified conjecture is proven for dimension $2$ and the trivial character (and for the sum over the non-trivial characters),  but also  for smooth quintics  of  dimension $3$ invariant under an involution.

Voisin's proof, apart from ingenious cycle theoretic constructions,   basically exploits the fact, as indicated above, that Chow groups of complete intersections  of the same degree but one dimension up become ``smaller''. The above remarks for complete intersections in weighted projective spaces motivate to look for suitable generalizations of her method. Although this can be done in any dimension, in this note I only propose  such a generalization for odd dimensional varieties (Theorem~\ref{result1}) since I only could find an interesting class of new examples for dimension $3$.  In this dimension  the result reads as follows (see Corollary~\ref{result2}):
\begin{theorem*}Let $W$ be a   smooth $(4+k)$-fold in  $\bP^{N+1+k}$, and let $Y\supset X\supset B$   smooth linear sections of dimensions $4,3, 2$ respectively. Suppose  that  the following assumptions hold:
\begin{enumerate}
\item $\chowhom_1(W)=0$
\item $\chowhom_0(Y)=\chowhom_{2}(Y)=0$; 
\item    there exists  a  family of  $1$-cycles on $X$ with $1$-dimensional parameter space, say $S$, such that the  Abel-Jacobi map  $S \to J^2(X)$ is onto ($J^2(X)$ is the intermediate Jacobian of $X$);
\item $h^{3,0}(X)=0=h^{1,0}(B)$ and $h^{1,1}(X)=1$.
 \end{enumerate}
Then, modulo  some further technical conditions (Assumption~\ref{assu}),  and provided $X$ is sufficiently general , $\chowhom_{0}(X)=0$. 
 \end{theorem*}
 Let me discuss the assumptions in the statement. 
 
 The third condition is well known to imply the generalized Hodge conjecture on $X$ (see \cite[Example 12.11]{PS}) and the condition $h^{3,0}(X)=0$ is  a consequence of $\chowhom_0(X)=0$ (by the generalized Mumford theorem).  
 Since it is often easier to prove that  $\chowhom_0(Y)=0$ for the variety $Y$ of which $X$ is a linear section it is natural to make this assumption. It will be explained below (Cor.~\ref{conseqfor4folds}) that in that case the Abel-Jacobi map for $2$-cycles cohomologous to zero  is injective, so  the vanishing of $\chowhom_2(Y)$ results if one knows for instance that $b_3(Y)=0$.   In this situation the Lefschetz theorem on linear sections then shows that $b_1(B)=0$ for the surface $B$.
 
One might wonder what the use of the variety $W$ is.  For this one goes back to condition (3). It remains true for the  entire  family of hypersurface sections of $Y$ with base $B$ and then gives a   correspondence, say $\beta$  from $Y$ to $B$ which one needs to show to vanish on the level of Chow-groups.  This is not immediately clear, but the functoriality of the construction shows that $\beta$ can be extended as a correspondence from $W$ to $B$ where $W$ is \emph{any}  variety of which $Y$ is a linear section. In that case $\beta$ factors over $W$ and by the philosophy of the linear section method $\chowhom_1(W)=0$ provided $k$ is large enough. 

For the moment  this method does not yet yield new examples of complete intersections, but there is one new example of a weighted complete intersection of general type. See  Example~\ref{subsec:wci}.
Nevertheless, the main theorem has potentially wider applicability, even for complete intersections and for that reason this case is explained in detail in \S~\ref{subsec:cis}. One can see this as a warming-up for the weighted complete intersection case in \S~\ref{subsec:wci}.
  
Finally a word about the technical conditions mentioned in the theorem. These  have to do with general position with respect to linear subspaces; one of these  is true  for complete intersections whose dimension is larger than half the dimension of the linear space (\S~\ref{sec:extend}) and -- provided suitably reformulated --  is  valid also for weighted complete intersections ; the second is true for all smooth projective varieties and this is the only statement which will be used for the example.

\subsection*{Acknowledgements} {\small I want to thank Robert Laterveer who posed several
inspirational questions which finally led to this paper and   Stefan M\"uller-Stach  and Vasudevas Srinivas who provided  technical assistance on various occasions.  I am also thankful to  the referee who spotted an  error  in the previous version and made me aware of  a couple of obscure statements and lines of reasoning.

\subsection*{Notation and conventions}  

\noindent In this paper   ``variety''  is  a reduced projective  scheme
over  the complex numbers so  it does not need to  be irreducible.

Let  $X$ be any variety. Then  $\chow_k(X)$ denotes  the Chow-group of $k$-cycles with $\Q$-coefficients modulo rational equivalence and $\chowhom_k(X)$, respectively $\chowalg_k(X)$, denotes   the subgroup of
$\chow_k(X)$ consisting of cycles homologous to $0$, respectively algebraically equivalent to $0$.
Occasionally one puts $\chow(X)=\oplus \chow_k(X)$.
One conventionally uses superscripts to denote codimension so
that $\chow^k(X)$ denotes the group of codimension $k$ cycles on $X$
modulo rational equivalence.

If a  variety is  considered as embedded in some fixed projective space, say $X\into \bP^N$ an important role is played by $\chowlin_k(X)$, the subgroup of $k$-cycles generated by $k$-dimensional linear sections and linear $k$-spaces contained in $X$. Of course $\chowlin_0(X)=\chow_0(X)$, but if $k>0$ the two may differ.

 For any variety $X$ the (singular) $m$-th cohomology group with \emph{rational} coefficients is denoted $H^m(X)$.  For a \emph{smooth} variety $X$ over $\C$   one denotes by
\begin{eqnarray*}
J^k(X)& = &H^{2k-1}(X;\C )/\left[F^kH^{2k-1}(X;\C )+H^{2k-1}(X;\Z)\right],\\
J^k(X)_\Q& = &H^{2k-1}(X;\C )/\left[F^kH^{2k-1}(X;\C )+H^{2k-1}(X)\right]
\end{eqnarray*}
the $k$-th intermediate jacobian, respectively the $k$-th intermediate jacobian   \emph{modulo torsion}. Here   $F^{\bullet}$ stands for the usual Hodge filtration. Given any $z\in \chowchom^k(X)$,
integration over a topological cycle whose boundary is $z$ defines
an element in the $k$-th intermediate jacobian (see \cite{Gr2} for
for details) and hence a homomorphism
$$
u^k_X: \chowchom^k(X) \to J^k(X)_\Q,
$$
the \emph{Abel-Jacobi map}.

\section{Correspondences} \label{sec:corresp}
 
\subsection{Generalities}

\noindent Let $X$ and $Y$ be smooth projective varieties. A \emph{correspondence
  from $X$ to $Y$} is a cycle, or equivalence class of cycles
on $X\times Y$. Write $\corr(X,Y)$ for the set of
$\Q$-correspondences from $X$ to $Y$ up to rational equivalence. If
$X,Y,Z$ are smooth, 
$\alpha\in\corr(X,Y)$, $\beta\in\corr(Y,Z)$, the composition
$\beta\comp\alpha\in\corr(X,Z)$ is defined as
$\beta\comp\alpha =(p_{13})_\ast[p_{12}^\ast \alpha\cdot
p_{13}^\ast \beta]$ where $p_{12},p_{13},$ and $p_{23}$ are the
projections of $X\times Y\times Z$ onto $X\times Y, X\times Z$, respectively
$Y\times Z$.

Any correspondence $\alpha$ from $X$ to $Y$ induces homomorphisms
between the Chow groups $\alpha_\ast:\chow(X)\to \chow(Y)$
by the formula $\alpha_\ast(u)=(p_2)_\ast[p_1^\ast u \cdot \alpha]$,
where $p_1$,   $p_2$ is the projection of $X\times Y$ to $X$
respectively  $Y$. This homomorphism does not necessarily preserve the degree.
If the components of $\alpha$ have the same dimension, say $\dim X+d$,
one says that $\alpha$ has degree $d$ since in this case
$\alpha_\ast$ induces homomorphisms
$$
\alpha_k: \chow_k(X)\to \chow_{k+d}(Y).
$$
Write $\corr_d(X,Y)$ for the degree $d$ correspondences from $X$
to $Y$.

The \emph{transpose} $\tr\alpha$   of $\alpha$ is obtained
by interchanging $X$ and $Y$. If $\alpha$ has degree $d$, $\tr\alpha$
has degree $\dim X-\dim Y+d$. In particular a degree $0$
correspondence between varieties of the same dimension induces \break
$\alpha_k:\chow_k(X)\to\chow_k(Y)$ and
$\tr\alpha_k:\chow_k(Y)\to\chow_k(X)$.

Any correspondence $\alpha$ between $X$ and $Y$ induces homorphisms
in homology by the same formula as the formula for
$\alpha_\ast$. In this formula $(p_1)_\ast$ is the usual map induced
in homology while $p_2^\ast$ is the Gysin map obtained from the map
in cohomology after applying Poincar{\'e} duality. If
$\alpha$ has degree $d$ the map is homogeneous of degree $2d$ and one
writes $\alpha_k : H_k(X)\to H_{k+2d}(Y)$.  This map is compatible with
the cycle class map $\cl_k: \chow_k\to H_{2k}(X)$. If one interchanges
the roles of $p_1$ and $p_2$ and passes  to cohomology one gets the
induced map in cohomology $\alpha^k:H^k(Y) \to H^{k-2d}(X)$. Under
Poincar{\'e} duality it coincides with $\tr\alpha_k$ which now is
  compatible with the cycle-class map $\cl^k:\chow^k(X) \to H^{2k}(X)$ for  cohomology.
It also induces a map between
intermediate jacobians $\alpha^k : J^k(Y) \to J^{k-d}(X)$ compatible with
the Abel-Jacobi maps, i.e.  there is a commutative diagram
\begin{equation}\label{eqn:FunforAJ}
\begin{matrix}
\chowchom^k(Y)
&\mapright{\alpha^k}&\chowchom^{k-d}(X)\cr
\!\!\!\!\!\!\mapdown{u_Y^k}&&\!\!\!\!\!\!\mapdown{u_X^{k-d}}\cr
J^k(Y)_\Q&\mapright{\alpha^k}&J^{k-d}(X)_\Q. \cr
\end{matrix}
\end{equation}
A  degree $d$ correspondence from $S$ to $X$  is nothing but a family  $Z_s$ of  $d$-cycles in $X$ parametrized by points $s\in S$. Fix $o\in S$. The $d$-cycles $Z_s-Z_o$ are homologous to zero and the map $S \to \chowhom_d(X)$ given by $s\mapsto [Z_s-Z_o]$ is the composition of  $S\to \chowhom_0(S)$ given by   $s\mapsto [s-o]$ and the   map $\chowhom_0(S)\to \chowhom_d(X)$ induced by the correspondence. Assembling the Abel-Jacobi images  of $ [Z_s-Z_o]$ then gives a morphism $S\to J^{n-d}(X)$, $d=\dim X$, which is also called the \emph{Abel-Jacobi map} associated to the correspondence.

\subsection{Degree zero correspondences between varieties of
the same dimension}

The following easy result is crucial for many of the results in this paper.
 \begin{lemma} \label{cycledecomposes} Let $\alpha$ be a degree zero
correspondence from
$X$ to $Y$, $\dim X=\dim Y=n$. Suppose that $\alpha$ can be
represented by a cycle having support on $V\times W$, $V$ a subvariety
of $X$, $W$ a subvariety of $Y$, $\dim V=v,\dim W=w$. Then 
\begin{enumerate}
\item  $\alpha_k:\chow_k(X)\to\chow_k(Y)$ is zero if $k<n-v$
or
$k>w$, 
\item  $\tr\alpha_k$ acts trivially on $\chow_k(Y)$ for
$k<n-w$ or $k>v$. 
\item  $\alpha_w$ acts trivially on $\chowhom_w(X)$ and
$\,\tr\alpha_v$ acts trivially on $\chowhom_v(Y)$.
\item  $\alpha_{w-1}|\chowhom_{w-1}(X)$ and
$\tr\alpha_{v-1}|\chowhom_{v-1}(Y)$ factor   over an abelian variety (mod torsion). Moreover  $\alpha_w$ as well as $\tr\alpha_{v-1}$ vanishes  on the respective Abel-Jacobi kernels, i.e. on $\ker(u_X^{n-(w-1)})$, $\ker(u_Y^{n-(v-1)})$ respectively.
\end{enumerate}
\end{lemma}
\proof
Let $\tilde V \to V$, $\tilde W\to W$ be
desingularizations and $i: \tilde V\to X$, $j: \tilde W\to Y$ be the
desingularization composed with the inclusions. Choose
$\tilde\alpha\in\corr(\tilde V\times\tilde W)$ mapping to $\alpha$.
With $p_1$ and $p_2$ the obvious projections, there is a commutative diagram
\begin{equation}
\begin{matrix}
\chow_{k+v+w-n}(\tilde V\times \tilde W)
&\rarrow{\tilde\alpha_{k+v+w-n}}{}{15mm} &\chow_k(\tilde V\times \tilde
W)\cr
\mapup{p_1^\ast}&&\mapdown{(p_2)_\ast}\cr
\chow_{k+v-n}(\tilde V)&&\chow_k(\tilde W)\cr
\mapup{i^k}&&\mapdown{j_k}\cr
\chow_k(X)&\rarrow{\alpha_k}{}{15mm}&\chow_k(Y).\cr
\end{matrix}\label{eqn:BasicDecomposition}
\end{equation}
Since $\chow_{k+v-n}(\tilde V)=0$ for $k<n-v$ and $\chow_k(\tilde
W)=0$ for $k>w$ the first assertion follows. The second follows from
the fact that $\tr\alpha_k$ factors over $\chow_{k+w-n}(\tilde W)$
and $\chow_k(\tilde V)$.

Since $\chowhom_v(\tilde V)=0=\chowhom_w(\tilde W)$ the third assertion
is clear.

Using the functoriality expressed by diagram~\eqref{eqn:FunforAJ} one deduces the fourth from the equalities  $\chowhom_{w-1}(\tilde W)=\pic^0(\tilde W)_\Q$, 
$\chowhom_{v-1}(\tilde V)=\pic^0(\tilde V)_\Q$  and the fact that the Abel-Jacobi kernels are zero for the Picard-varieties. 
\endproof
 
To apply this Lemma, one needs a very strong decomposition  for the correspondences.  If the varieties have small Chow groups this can in some cases be achieved using  the  Bloch-Srinivas method \cite{B-S}. For this note  the following variant is useful.   
\begin{prop} \label{factor} Let $X,Y$ be varieties of the same dimension $n$ and let $\alpha\in \corr_0(X,Y)$ be a degree $0$ correspondence. Suppose that the image of $\alpha_0:\chow_0(X) \to \chow_0(Y)$ is supported on a subvariety $W\subset Y$. Then in $\chow_n(X\times Y )$ one has a decomposition $\alpha=\alpha^{(1)}+\alpha^{(2)}$ with $\alpha^{(1)}$ supported on $X\times W$ and $\alpha^{(2)}$ supported on $D\times Y$ where $D$ is some (possibly reducible) divisor on $X$.
\end{prop}
This variant can be proved by copying the proof of  \cite[Prop. 1]{B-S} which is the case where $\alpha$ is the diagonal inside $X\times X$.   See also \cite[Corollary 10.20]{V2}. 
A first consequence of Prop.~\ref{factor}, but in fact of \cite[Prop. 1]{B-S} is the following result:
 \begin{cor}
\label{diag} Let   $X$ be a smooth
variety of dimension $n\ge 2$ with $\chowhom_0(X)=0$. Then there
exists a smooth variety $V$ of dimension $n-2$ and $\beta
\in\chow_{n-1}(X\times V)$ such that $\tr\beta_k\comp\beta_k=\id$
for  $1\le k\le n-1$.
\end{cor}
\proof
By Prop.~\ref{factor}  (or by  \cite[Prop. 1]{B-S}) there is  a decomposition
of the diagonal   inside $\chow(X\times X)$ as  $\Delta= x\times X
+\gamma$, $x\in X$ and such that $\supp(\gamma)\subset X\times D$, $D\subset X$
a divisor, which one may assume to be very ample.

If $k<n$ the first factor acts trivially on $\chow_k(X)$, an then  $\gamma$ acts as the identity. The support of $\gamma$ cannot be
contained in $X\times E$, with  $E$ a variety  of dimension $\le n-2$: the
identity factors over $\chow_k(\tilde E)$, with $\tilde E$ a
desingularization of $E$ and so is trivial for $k\le n$, while
$\chow_{n-1}(X) \not=0$. In particular, taking a hyperplane section
$D\cap H$ of $D$, the cycle
$$
\beta=\gamma\cdot X\times H
$$ 
has support of dimension $n-1$. Choose a desingularization $\tilde
D\to D$ which contains a desingularization $V\to D\cap H$ of $D\cap
H$. Let $i:\tilde D\to X$ and $j:V\to\tilde D$ be induced by the
obvious inclusions and let $\tilde\gamma\in\corr(X,\tilde D)$ be a
cycle class mapping to $\gamma$. By assumption,
$i_\ast\comp\tilde\gamma=\gamma$ acts as the identity.  The same argument for
the transposed cycles proves that$\tr\tilde\gamma\comp i^\ast=\tr\gamma $ acts
as the identity  provided $k>0$. One may assume that $H\in|D|$ so that $i_*\comp i^*$ and
$j^*\comp j_*$ both equal cup product with the first Chern class of
$c_1(D)=c_1(H)$, one arrives at a commutative diagram
 
\vspace*{2cm}
\hspace{5mm}
\setlength{\unitlength}{1mm}
\begin{picture}(0,0)
\put(2,13){\vector(4,-3){10}}
\put(87,5){\vector(4,3){10}}
\put(37,4.5){\line(1,0){24}}
\put(37,3.5){\line(1,0){24}}
\put(9,9){$_{\beta_k}$}
\put(94,9){$_{\tr\beta_k}$}
\put(-14,9){
$\begin{array}{cccccccc}
\chow_k(X)& \rarrow{}{\tilde\gamma}{5mm} &\chow_k(\tilde D)
&\rarrow{}{i_\ast}{5mm} &\chow_k(X)\rarrow{}{i^\ast}{5mm}  
&\chow_{k-1}(\tilde D)&\rarrow{}{\tr\tilde\gamma}{5mm}&\chow_k(X)\\ 
&&\darrow{j^\ast}{}{4mm}&&&\uarrow{}{j_\ast}{4mm}\\
\noalign{\vskip-1mm} 
&&\chow_{k-1}(V)&&& \chow_{k-1}(V) \\
\end{array}$}
\end{picture}

\vspace{3mm}

\noindent It shows that $\tr\beta_k\comp\beta_k$ is the identity.
\endproof
I need another consequence of  \cite[Prop. 1]{B-S}:
\begin{cor}[ \protect{
\cite[theorem 1. p. 1238]{B-S}}]
 \label{conseqfor4folds}
Suppose that $\chowhom_0(X)=0$. Then  the    Abel-Jacobi map   
$$u_X^{2}: \chowchom^{2}(X) \to J^{2}(X)_\Q$$ 
is injective; in particular, if  $H^3(X)=0$, one has $\chowchom^2(X)=0$. 
\end{cor}

\subsection{Reflexive vanishing for correspondences}

\begin{definition}\label{reflex}
Let $X$ be a smooth variety.  It has the \emph{$m$-th reflexive vanishing property with respect to correspondences} if the following property holds:
\[
\leqno{\rm C}(X)_m
\left\{ \begin{array}{ll}
\text{For any smooth $V$ with $\dim V=\dim X$ and any degree $0$ } \\
\text{correspondence  $\alpha$ from $V$ to $X$,   if the induced map}\\
\text{$\alpha_m:\chowhom_{m}(V)\to \chowhom_{m}(X)$ vanishes  then  $\tr\alpha_m=0$. }
\end{array}\right.
\]
 \end{definition}

This  reflexive vanishing} property is difficult to test but  for $m=0$ one can use the Bloch-Srinivas method from  the previous section.
\begin{prop} \label{reflexforsurfs} Let $X$ be a surface with $\alb(X)=0$. Then $X$ has the $0$-th reflexive vanishing property with respect to correspondences.
\end{prop}
\proof
Suppose that $\alpha\in\corr_0(V,X)$ has the property that $\alpha_0=0$ on $\chowhom_0(V)$. This means that the image of $\alpha_0$ is supported on   points. By Prop.~\ref{factor} one can write $\alpha=\alpha^{(1)}+\alpha^{(2)}$ with $\alpha^{(1)}$ supported on $V\times$ [points], while $\alpha^{(2)}$ is supported on [divisor]$\times X$. Now apply Lemma~\ref{cycledecomposes} to see that $\tr\alpha_0$ factors over an  abelian variety and hence over $\alb(X)=0$.
\endproof

\begin{rmq} The referee remarks that the  property {\rm C}$(X )_0$  is satisfied if and only if for any variety 
$Y$ of dimension $< \dim X$, and any correspondence $\alpha\in\corr^0(X, Y )$, $\alpha_0: 
\chowhom_0(X ) \to \chowhom_0(Y)$ is trivial. Indeed, one implication follows immediately from the  preceding proof. Conversely,  if C$_0(X)$ holds, pick a smooth projective variety $E$ such that $\dim E+\dim Y =\dim X=n$; Lemma~\ref{cycledecomposes} shows that if $\alpha\in \corr_0(X,Y)$ is viewed as a correspondence $\beta$ from $X$ to $Y\times E$ (by choosng a point $e$ of $E$ and identifying $Y$ with $Y\times e$), one has $\tr \beta_0=0$ and hence, by reflexivity, $\beta_0=0$, whence $\alpha_0=0$.
\\
In terms of the ``na\"ive filtration'' on $\chow_0 (X)$  introduced by Voisin in \cite{V3}  this means by definition that 
\[
F_{\text{\rm  na\"\i ve}}^i \chow_0 (X ) =F_{\text{\rm  na\"\i ve}}^n \chow_0 (X ), 0 < i 
\le  n, n = \dim X.
\]
\end{rmq}

\section{Varieties with linearly generated Chow groups} \label{sec:linchow}

\subsection{Complete intersections in (weighted) projective space} 

\begin{definition}\rm  \label{islingen} Let $X\subset\bP^N$. 
Let me say that $\chow_k(X)$ is \emph{linearly generated} if it is  generated by
linear sections and linear subvarieties of $\bP^N$ contained in $X$.
\end{definition}

\begin{prop}[\cite{E-L-V}]   \label{ComplInts} Let  $X$ be a complete intersection of
multidegree
$(d_1,..\allowbreak.., d_r)$ in $\bP^{n+r}$. Assume that $d_1\ge d_2\dots \ge
d_r\ge 2$. For any $s$ such that
\begin{equation}
\sum {d_j+s\choose s+1} \le n+r-s, \label{Threshold}
\end{equation}
the Chow groups $\chowhom_k(X)$ vanish for $k\le s$, while
$\chow_{s+1}(X)$ is linearly generated: $\chow_{s+1}(X)=\chowlin_{s+1}(X)$.

If $d_1\ge 3$ or $r\ge k+1$ the vanishing result remains true  if one replaces the
right hand side of (\ref{Threshold}) by $(n+r)$; for an intersection of $r\le
k$ quadrics one can replace the right hand side by $n+2r-k-1$. 
\end{prop}

From these examples one finds new ones in weighted projective spaces as follows.
Recall that a weighted projective space $\bP(a_0,\dots,a_n)$ is the
quotient of $\C^{n+1}\setminus  \set{0}$ by the action of $\C^\ast$ given by
$t\cdot(z_0,\dots,z_n)=(t^{a_0}z_0,\dots,t^{a_n}z_n)$.  Equivalently, letting
$\zeta_i$ be a primitive $a_i$-th root of unity,  let the generator of the
$i$-th factor of  
$G=\Z/\Z_{a_0}\times\cdots\times \Z/\Z_{a_n}$ act on $\bP^n$ by multiplying the
$i$-th coordinate with $\zeta_i$, the quotient is $\bP(a_0,\dots,a_n)$.

One has two different concepts of   a linear space inside a weighted projective space
$P=\bP(a_0,\dots,a_n)$. The first kind is the   image  of  any linear space 
in $\bP^n$ and is called \emph{weighted linear subspace}.  Note that a weighted codimension $1$ subspace  need not be given by an equation in homogeneous coordinates. Those for which an equation can be found  will be called \emph{weighted hyperplanes}.

The second kind, called \emph{generalized linear spaces}, come from the following construction. Set $M+1=\sum a_i$
and introduce homogeneous coordinates $X_i^j$, $i=0,\dots,n$, $j=1,\dots,a_i$ 
on $\bP^M$. For any choice $J=(j_0,\dots,j_n)$, $1\le j_k \le a_k$,
define $L_J$ as the linear codimension $n+1$ subvariety of $\bP^M$ given
by  the equations $\xi_k^{j_k}=0,k=0,\dots,n$ and set $Q=\bP^M\setminus
\bigcup_J L_J$. Define  $\varphi: Q\to P$ by  $\varphi(\dots,X_i^j,\dots)=(\dots,X_i^1\cdots
X_i^{a_i},\dots)$. If one embeds $\bP^n$ with homogeneous coordinates
$(x_0,\dots,x_n)$ in $Q$ by setting $X_i^j=x_i$, the restriction of $\varphi$
to $\bP^n$ is just the natural quotient map $\bP^n\to P$. It is clear that $Q$
only contains complete linear spaces of dimension $\le n$.  
 Any linear $k$-space inside $Q$ meeting the general fibre  of $\varphi$ in at
most points maps to a variety in $P$ of    dimension $k$ and is called a
\emph{generalized $k$-plane} of $P$.

Since $z_i$ has weight $a_i$,  a homogeneous polynomial in the $z_i$ of weighted degree $d$ defines a hypersurface in $P=\bP(a_0,\dots,a_n)$. It is called \emph{quasi-smooth} if the corresponding
hypersurface in $\C^{n+1}$ has $0$ as its only singularity. A similar
definition holds for complete intersections. 
  The construction of generalized linear spaces  can be used to prove:
\begin{prop}[\cite{Le}] \label{lowchow} Let  $X$ be a complete intersection of multidegree  \allowbreak
$(d_1,\dots,d_r)$ in $P=\bP(a_0,\dots,a_{n+r})$. Put $M+1=\sum_{j=0}^{n+r} a_j$.
Assume that $d_1\ge d_2\dots \ge d_r\ge 2$. For any $s$ such that
$$
\sum {d_j+s\choose s+1}\le M-s,
$$
the Chow groups $\chowhom_k(X)$ vanish for $k\le s$, while
$\chow_{s+1}(X)$ is generated by generalized linear subspaces contained in $X$.
\end{prop}
\begin{example}\rm  I am particularly interested in the following weighted projective space $\bP(2^a,3^b)$, where the notation means that one takes $a$ weights to be equal to $2$ and $b$ weights equal to $3$.
Hence $M=2a+3b-1$. The complete intersections $X$  given by homogeneous equations of degree $6$ are the simplest ones.  For these  $\chowhom_0=0$ whenever $r\le \frac 1 6 (2a+3b-1)$ and $\chowhom_1=0$ whenever $r\le \frac{1}{21} (2a+3b-2)$.  Rephrasing this in terms of the dimension $n$ of the intersection, the first condition becomes
\[
4a+3b\le 6  (n+1).
\]
One sees that for $3$-folds this  gives $4a+3b\le 24$  and for $4$-folds  $4a+3b \le 30$. The special case $a=4$, $b=3$ is treated later. For the  $3$-fold  $X=X_{6,6,6}$ this result does not give $\chowhom_0(X)=0$ but for the fourfold $Y=Y_{6,6}$ it does give $\chowhom_0(Y)=0$.
\end{example}

\subsection{Extending linear cycles in pencils}\label{sec:extend}

\begin{lemma} \label{cutout} Let $P=\bP^{n+r}\subset P'=\bP^{n+r+1}$,
$L\subset P'$ a linear space   of dimension $k+1$, and let $M=L\cap P$  be
the corresponding linear space of dimension $k$ in $P$. Let $X\subset P$ be an $n$-dimensional smooth complete intersection
containing $M$.  
\newline
{\rm (i) } Suppose that $ k
<\frac 12 n $; then   there exist  a smooth complete intersection
$Y$ in $P'$ of the same multidegree as $X$, which contains   $L$  and for which
$Y\cap P=X$.
\newline {\rm (ii)} If $k <\frac 12 (n-1)$ the general hyperplane  of  $P'$ passing through $L$ cuts $Y$ in a smooth hypersurface.
\end{lemma}
\proof
(i) Put $N=n+r$. and choose
homogeneous coordinates
\[
\set{X_0,\dots,X_N,X_{N+1}}
\] 
in $P'$ so that $P$ is given by $X_{N+1}=0$. Let $X$ be given by the equations  $F_i(X_0,\dots,X_N)=0$, $i=1,\dots,r$. Suppose that $\deg F_i=d_i$.  Let $G_i=0$ be a hypersurface of degree $d_i-1$ in $P'$ containing $L$ and consider the complete intersections   inside $P'$ given by equations $F_i+X_{N+1}G_i=0$, $i=1,\dots,r$. They all contain  $L$  which is in fact  the base locus of the linear system of these complete intersections. The generic member of this system can only be singular at points on $L$. I will show that   specific members cannot have singularities at all along $L$ which therefore implies that for generic choices of $G_i$ the resulting complete intersection satisfies the demands. 

To find such a complete intersection,
pick any codimension $r$ linear space of $P$, say given by the linear equations $H_1=\cdots=H_r=0$ and put $G_i=H_iK$ where $K(X_0,\dots,X_{N+1})$ is of degree $d_i-2$ and $K=0$ contains $L$. Suppose that $q=(p_0,\dots,p_N,p_{N+1}) \in L$ is a singularity on the corresponding complete intersection.
Let $p=(q_1,\dots,q_N,0)\in M$. Then $p\in X$ (since  $X$ contains $M$). Let $\nabla$ be the gradient with respect to $(X_0,\dots,X_N)$.  Then the vectors
\[
\nabla_pF_i+q_{N+1}K(q) \nabla_pH_i, \quad i=1,\dots, r
\]
have to be dependent, say $-q_{N+1} K(q)  \nabla_p \sum_ i \lambda_i  H_i = \nabla_p  \sum \lambda_i F_i$. This implies that the hyperplane $ \sum_ i \lambda_i H_i =0$ contains $T_pX$ since the latter is the common annihilator of  $\nabla_pF_1,\cdots,\nabla_pF_r$ while any hyperplane given by an equation $H=0$ is the annihilator of $\nabla_p H$. So I have found a hyperplane of $P$ which contains   $M$ and $T_pX$. 
It suffices to find a hyperplane for which this is not possible: 
\begin{claim} There exists a hyperplane $H\subset P$ passing  through $M$  which does not contain any of the tangent planes to $X$ at points $p\in M$.
\end{claim}
To prove this, fix   a linear $(N-k-1)$-dimensional subspace $M'$   of $P$ disjoint from $M$. 
Fix a hyperplane $H$ containing  $M$.  Such a hyperplane varies in the projective space $ (M')^*$ of hyperplanes in $M'$. 
 Let $I=\sett{(p,H)\in M\times  (M')^*}{ T_pX\subset H}$.  
 The variety $I$ has dimension  $ k+(r-1)$ since the projection $I\to M$ is surjective with fibres of dimension  $(r-1)$ since $T_pX\cap M'$ has codimension $r$ so the hyperplanes containing it form a projective space of dimension $(r-1)$.  It suffices to show that  the projection $I\to (M')^*$ is not surjective. For that it suffices to observe  that  $\dim I= k+(r-1)< \dim (M')^*= N-k-1=n+r-k-1$ since $2k< n$. 
 
\medskip
\noindent (ii) The argument is similar: replace $M$ by $L$, $M'$ by a linear subspace  $L'\subset P'$ disjoint from  $L$ with $\dim L'= N-k-1$.The  hyperplanes $H$ containing  $L$ vary in the projective space $ (L')^*$ of hyperplanes in $L'$ and if $H\cap Y$ is singular at $q\in L$ one has $T_q(Y)\subset H$. Consider $J=\sett{(q,H)\in L\times  (L')^*}{ T_qY\subset H}$. Since $\dim J= k+1+(r-1)=k+r$ and $\dim J= k+r < n+r-k-1= \dim (L')^*$ the projection $I'\to (L')^*$ is not onto and if $H$ is a  hyperplane corresponding to a   point $[H]$  not in the image, it cannot contain any tangent plane of $Y$ and the intersection is smooth.
\endproof

\begin{rmk} \rm 1) \label{cutout2} The first assertion in Lemma~\ref{cutout} can be generalized to quasi-smooth weighted complete intersections with respect to weighted linear subspaces. 
\\
2) The proof of the second assertion shows that it has nothing to do with complete intersections: it is valid for any smooth $Y\subset P'$ containing a linear space $L$ which satisfies $\dim L <\frac12 \dim Y$.
\end{rmk}

\section{Normal functions of pencils for which the generalized Hodge conjecture holds}
 Consider the middle cohomology groups for a smooth  $(2m+1)$-fold $X$. By \cite[Example 12.11]{PS}  the  generalized Hodge conjecture is guaranteed by the   following condition on  $m$-cycles:
\[  \left\{ \begin{array}{ll} \text{There exist  $S$ and  $\alpha\in \corr_m(S,X)$ such that its}\\
\text{Abel-Jacobi map $\alb(S) \to J^{m+1}(X)$  is onto}.
\end{array}
\right.
 \leqno{\rm GHC}_m(X) 
 \]  
It is implicitly assumed that the relative cycle defining $\alpha$ restricts  on the fibres $\set{s}\times X$ to cycles which are homologous to zero. By Lefschetz theorem on the hyperplane section and functoriality, one can always suppose that  $S=C$, a smooth  curve and then  $\alpha\in\chow_{m+1}(C\times
X)$ and the induced map $\alpha^1: J^1(C)_\Q\to J^{m+1}(X)_\Q$ is onto. 
 It follows that  $\alpha^1\comp\tr\alpha^1$ is a non-zero endomorphism of $J^{m+1}(X)_\Q$.  Suppose that  
 \begin{equation}\label{eqn:mult}  \End(J^{m+1}(X)_\Q)= \Q\cdot \id,
\end{equation} 
then one has 
\begin{equation} \label{eqn:mult2}
\alpha^1\comp\tr\alpha^1=r(\alpha)\cdot \id,\quad r(\alpha)\not=0.
\end{equation}
 
The main result is:
\begin{prop} \label{secondobs} Let $Y\subset \bP^{N+1}$ be a smooth $(2m+2)$-fold  with 
\begin{equation*}  
\chowhom_{m+1}(Y)=0
\end{equation*}
and let   $j:X\into Y$  be a smooth hyperplane section.   Suppose that $H^{2m+1}(X)$ has Hodge level $\le 1$, i.e.  
\begin{equation*}  H^{2m+1}(X)=H^{m,m+1}(X)+H^{m+1,m}(X),
\end{equation*}
and that, moreover,  \eqref{eqn:mult} holds  as well as    assumption {\rm GHC}$_m(X)$.
Then there exists    $r\in\Q^\ast $ such that   for any  $d' \in \chowhom_m(X)\cap j^* \chow_{m+1}(Y)$ one has $\alpha\comp \tr\alpha (d')= rd'$ \footnote{A priori $r$ depends on $d'$, but  the proof shows that this is in fact  not the case.}.
\end{prop}
To prove this Proposition, one lets $X$ vary in a hyperplane pencil. So let $L_{N-1} \subset \bP^{N+1}$ a linear subspace of codimension $2$.  The linear pencil of hyperplanes passing through $L_{N-1}$ is parametrized by a line $L\subset \bP^{N+1}$ disjoint from $L_{N-1}$. The corresponding hyperplane sections are denoted by $X_t=Y\cap H_t$, $t\in L$. 
Set    $\tilde Y=\sett{(y,t)\in Y\times L}{y\in X_t}$ and let $\tilde f:\tilde Y \to L$ be induced by the second projection. Let  $U\subset L$ be the Zariski-open subset such that $X_t$ is smooth for $t\in U$.  Set $Y':= \tilde f^{-1}U$, $f=\tilde f |Y': Y'\to U$. 
Start with any  $(m+1)$-cycle $Z$     on $Y$ whose intersection $Z_t$ with a  smooth fibre  $X_t$  is   homologous to zero on that fibre.  Then $t\mapsto [Z_t]\in J^{n-m}(X_t)$ defines the \emph{normal function} $\nu_Z$ associated to $Z$ with cohomology invariant $\delta(\nu_Z)\in H^1(U,R^{2m+1}f_\ast\Q)$. In the case under consideration, the Leray spectral sequence for $f$ gives an
isomorphism
\begin{equation}
\tau: \ker\bigl(H^{2m+2}(Y')\to
H^0(U,R^{2m+1}f_\ast\Q)\bigr)\mapright{\cong} H^1(U,R^{2m+1}f_\ast\Q)
\label{two}
\end{equation}
and it is well known   \cite[(3.9)]{Zu} that $\tau$ maps $[Z]_{Y'}$, the
cohomology class of $Z$ on $Y$, which by assumption is in the left hand side,  to
the cohomology invariant $\delta(\nu_Z)$. So   $[Z]_{Y'}=0$ is zero precisely when  $\delta(\nu_Z)=0$.

The starting observation is:
\begin{lemma}  \label{startobs} In the above setting, suppose that $\delta(\nu_Z)=0$. Then for all $t\in U$ one has $Z|X_t=0$ in $\chowhom_{m}(X_t)$.
\end{lemma}
\proof
The cycle  $Z$  pulls back to  $\tilde{Y}$ and
restrict to a cycle  on $Y'$ which is denoted by the same symbol. Since
$\tau$ (see  \eqref{two}) is an isomorphism, its  cohomology class $[Z]_{Y'}$
 as a class on $Y'$  vanishes.  Let $\Sigma$ be the union of the singular fibres so that $Y'=\tilde Y\setminus\Sigma$ and let $k:Y'\into \tilde Y$ be the inclusion. There is a commutative diagram
\[
\begin{array}{ccccccc}
\ker(k_*) &\to &\chow_{m+1}(\tilde
Y) &\mapright{k_*} &\chow_{m+1}(Y') &\to &0\\
\mapdown{\cl_1}&&\mapdown{\cl_2}&&\mapdown{\cl_3}\\
\bigl(H^{2m+2}_\Sigma(\tilde Y)\bigr)_\alg&\to& H^{2m+2}_\alg(\tilde Y)&\to&H^{2m+2}(Y',\Q).  
\end{array}
\]
The subscript ``alg'' stands for the subspace generated by classes of algebraic cycles.  In particular,  $\cl_1$ surjective, and the assumption $\chowhom_{m+1}(Y)=0$  implies that $\cl_2$ is an isomorphism on the subgroup of cycles coming from $Y$.
The diagram then shows that  $[Z] \in\chow_{m+1}(\tilde Y)$ is
supported on the singular fibres. The restriction of $[Z]$ to a smooth
fibre therefore  vanishes.  
\endproof

\noindent \textit{Proof of Proposition~\ref{secondobs}.}
Take an $(m+1)$-cycle $Z'$  on $Y$ which restricts  to a cycle on $X$ representing $d'$. Assume that  $Z'$ is chosen so that it is  transversally intersected by any smooth member $X_t$ of the pencil  in an $m$-cycle   ${Z'}_t$  homologous to zero on $X_t$.  The assumption GHC$_m(X_t)$ provides a family of curves $C_t$, $t\in U$ and  correspondences  $\alpha_t$ on $C_t\times X_t$. Possibly after some base change $U'\to U$ these fit into a global relative family $C_{U'}\to U'$ and a relative correspondence $\alpha_{U'}$  from $C_{U'}$  to $Y'/U'$. It induces a global endomorphism of the family of intermediate Jacobians $J^m(X_t)$ which must be a homothety since this is the case  for generic $t$. Hence the   numbers  $r(\alpha_t)\in\Q$ from \eqref{eqn:mult2}  for those $t\in U$ for which  \eqref{eqn:mult} holds,   are independent of $t$, say $r(\alpha_t)=r(\alpha)$. The Abel-Jacobi maps for the family $Z'_t$ yield the normal function $\nu_{Z'}$ while the normal function for the $(m+1)$-cycle $\alpha\comp \tr \alpha(Z')$ is just $t\mapsto  \alpha^1\comp\tr\alpha^1(\nu_{Z'}(t))$ which  is equal to $r(\alpha)\nu_{Z'}(t)$.   Lemma~\ref{startobs} applied to the difference  $Z= \alpha\comp \tr \alpha(Z')-r(\alpha)\cdot Z' $ then completes the proof. \qed

\section{Controlling  Chow groups using  pencils of hyperplane sections  }

The aim is to show the following result:
\begin{theorem}  \label{result1} Let $W\subset  \bP^{N+1+k}$ be a fixed $(2m+2+k)$-fold.  Suppose moreover, that
\begin{enumerate}
\item $\chowhom_m (W)=0$ ; 
\item $\chowhom_0(Y)=\chowhom_{m+1}(Y)=0$ for every smooth codimension $k$  linear section $Y$ of $W$; 
\item for any  smooth hyperplane section $X$ of $Y$  one has 
$h^{m-1,m-1}(X)=h^{m,m}(X)=1$,   $H^{2m+1}(X)$ has  Hodge level $\le 1$, and, moreover   condition  {\rm GHC}$_m(X)$ holds;
\item    any smooth  linear section $B$ of $Y$ of codimension $2$ has the $(m-1)$-st reflexive vanishing property with respect to correspondences;
\end{enumerate}
Then, provided $X$ is sufficiently general, and the technical condition below (Assumption~\ref{assu}) holds, one has $ \chowlin_{m-1}(X)\cap \chowhom_{m-1}(X)=0$. 
\end{theorem}
For $m=1$ this can be simplified   since  Prop.~\ref{reflexforsurfs}  states that for $0$-cycles on surfaces with $q=0$ the reflexive vanishing holds which yields:
\begin{cor}  \label{result2} Let $W$ be a   smooth $(4+k)$-fold in  $\bP^{N+1+k}$, and let $Y\supset X\supset B$   smooth linear sections of dimensions $4,3, 2$ respectively. Suppose that Assumption~\ref{assu} holds. Suppose moreover, that  
\begin{enumerate}
\item $\chowhom_1(W))=0$
\item $\chowhom_0(Y)=\chowhom_{2}(Y)=0$; 
\item     {\rm GHC}$_1(X)$ holds;
\item $h^{3,0}(X)=0=h^{1,0}(B)$  and $h^{1,1}(X)=1$.
 \end{enumerate}
 Then, provided $X$ is sufficiently general, one has  $\chowhom_{0}(X)=0$. 
\end{cor}

The technical assumption alluded to above is inspired by Lemma~\ref{cutout}: 
\begin{ass} \label{assu} Let $s=m-1$ or $s=m$. Let $M \subset X$ be a given linear space of dimension $s$. Then\\
1) there exists a linear space $L=L_{s+1}\subset \bP^{N+1+k}$ of dimension $(s+1)$ and a smooth codimension $k$ linear section $Y $  of $W$ containing $L$ such that  $X$ is a hyperplane section of $Y$ with  $X\cap L=M$;   \\
2) the generic  hyperplane section $X' $  of $Y$  containing  $L$  is smooth.
\end{ass}

\noindent\textit{Proof of Theorem~\ref{result1}.}
Consider a class  of the form $c=\deg(X)\cdot [M]- [H^{m+2}\cap X]$,  with $M$ a linear space of dimension $(m-1)$ contained in $X$ and $H$ a hyperplane. By  hypothesis  3) of the theorem such $c$ belong  to  $ \chowhom_{m-1}(X)$ and these classes generate $\chowlin_{m-1}(X)\cap \chowhom_{m-1}(X)$. So  it suffices to prove that $c=0$.  Apply assumption~\ref{assu} as follows. Starting   with the  linear space $M \subset X$ one finds $Y$ containing $L$  as in 1,  and a smooth $X'$ as in 2). Next,  one applies  assumption 1) to $(X',L)$ to find $(Y',L')$.
The situation can be summarized in  the diagram
\[
\begin{matrix}
&& X' & \mapright{j'}&  Y' \\
&& \mapdown{k}\\
X& \mapright{j}&Y,
\end{matrix}
\]
where all arrows are inclusions of hyperplane sections $Y,Y'\subset W$,  a   $(2m+2+k)$-fold $W\subset  \bP^{N+1+k}$. By construction one has $j^*[L]=[M]$.  Put $d= \deg(Y)\cdot [L]- [H^{m+2}\cap Y]$. By hypothesis  3) of the theorem  and Lefschetz' hyperplane theorem $h^{m,m}(Y)=h^{m,m}(X)=1$ and hence $d \in \chowhom_m(Y)$. Clearly $j^*d=c$. The smooth hyperplane section $X'$ of $Y$ contains  $L$ and hence $L$ defines a class   $[L]\in \chow_m(Y')$ supported on $X'$ and so     $d=k_*[d']$ for some $d'\in\chowlin_m(X')\cap \chowhom_m(X')$ and since $L'\cap X'=L$,
for $e'=\deg(Y') [L']-  [H^{m+1}\cap Y' ] \in  \chow_{m+1}(Y')$ one has $(j')^*e'=d'$.\footnote{Since  $h^{m+1,m+1}(Y')\not=1$ in general, there is no a priori reason why $e'$  should be  homologous to zero; one can only say that the class is primitive.}
The main result of the previous section, Prop.~\ref{secondobs}   can then be applied to the class $d'\in \chowhom_m(X')$ and to the inclusion $j':X'\into Y'$ (instead of $j$).   
Hence  
\begin{equation*}
\alpha\comp\tr\alpha (d')=r d' . 
\end{equation*}
The varieties $X$ and $X'$ fit in a pencil of hyperplanes of $\bP^{N+1}$ passing   through a certain linear subspace $L_{N-1}\subset \bP^{N+1}$ of codimension $2$ and  parametrized by a line $\bP^1\subset \bP^{N+1}$ disjoint from $L_{N-1}$. The base locus $B=Y\cap L_{N-1}$ is smooth and   $\tilde Y=\sett{(y,t)\in Y\times \bP^1}{y\in X_t}$ is the blow-up of $Y$ in $B$  with blow-down map  $\sigma: \tilde Y \to Y$  induced by the second projection.  The first projection $p:\tilde Y \to \bP^1$ is the pencil. 
It gives a natural embedding of $X$ into $\tilde Y$ as  some  fibre  of  $p$.  For the smooth members $X_t, t\in U$ of the pencil the property GHC$_m(X_t)$ holds.  Hence  there is a relative correspondence $\alpha_U\in \corr(C_U,\tilde f^{-1}(U))$ extending to  $\tilde \alpha \in \corr(C_{\bP^1} ,\tilde Y)$.    Let $q:C_{\bP^1}\to \bP^1$ be the natural projection.
By construction 
\begin{equation*} 
\tilde\alpha\comp\tr\tilde\alpha| X_{p(y)} = \alpha_{q(y)} \comp\tr\alpha_{p(y)}.
\end{equation*}
Denote    the proper transform in $\tilde Y$ of a class  of a cycle  on $Y$  by placing a tilde over the class. The above equation  then shows:
  \begin{equation}
\tilde\alpha\comp\tr\tilde\alpha   (\tilde d)=r tilde d, \quad d=k_*d' .  \label{eqn:transposing2}
\end{equation} 
Composing $\tilde\alpha$  with the correspondence
from $\sigma: \tilde Y \to Y$ given by the blow down morphism one obtains $\tilde\alpha^{(1)}=\sigma_*\comp \tilde\alpha\in \corr_m(C_{\bP^1},Y)$. Composing it with the
correspondence from $\tilde Y$ to $B$ defined by inclusion of the exceptional divisor $E\into  \tilde Y$ followed by $\sigma|E$, one gets the   correspondence $\tilde\alpha^{(2)}\in \corr_{m-1}(C_{\bP^1},B)$. The induced homomorphisms
$$
\begin{array}{ll}
\tilde\alpha^{(1)}&:\chow_0(C_{\bP^1})\to\chow_m(Y)\cr
\tilde\alpha^{(2)}&:\chow_0(C_{\bP^1})\to\chow_{m-1}(B) 
\end{array}
$$
come from the decomposition 
\begin{equation} \label{eqn:decblow} 
\chow_m(\tilde
Y)\cong\chow_m(Y)\oplus\chow_{m-1}(B).
\end{equation}
On the other hand, any class in $\chow_m(\tilde Y)$  such as $\tilde d$ can be decomposed according to the decomposition \eqref{eqn:decblow}:  $\tilde d=d  + s$ and \eqref{eqn:transposing2} implies
\[
r j^\ast d  =
j^\ast\comp\tilde\alpha^{(1)}\comp\tr\tilde\alpha^{(1)}d  +
j^\ast\comp\tilde\alpha^{1}\comp\tr\tilde\alpha^{(2)} s.
\]
Put  $u= \tr\tilde\alpha^{(1)} d\in\chow_0(C_{\bP^1})$. Since $d$ is supported on $X'$, also  $d':=\tilde\alpha(u)\in\chow_m(\tilde Y)$ is supported on $X'$ so that $j^*\sigma_*d'$  is supported on $X'\cap X=B$. In other words -- with $i:B\into X$ the inclusion --- one can write 
\begin{equation}\label{eqn:transfer}
j^*\sigma_*d'=i_*s', \quad r'\in \chow_{m-1}(B).
\end{equation}
This means that  
\[
d' =\sigma_*d'+ s'
\]
is the decomposition of $d'$ according to \eqref{eqn:decblow}. Hence \eqref{eqn:transfer} translates into
\[
j^*\tilde\alpha^{(1)}(u)=i_*\tilde\alpha^{(2)}(u),\quad u=\tr\tilde\alpha^{(1)} d
\]
and hence
\[
rj^\ast d  = i_\ast\comp\tilde\alpha^{(2)}\comp\tr\tilde\alpha^{(1)}d +
j^\ast\comp\tilde\alpha^{1}\comp\tr\tilde\alpha^{(2)} s.
\]
One deduces:
\begin{observation}   Let $d\in \chow_{m }(Y)$.  Then the $(m-1)$-cycle $j^\ast d$ on $X$ is rationally equivalent to zero  provided 
\begin{eqnarray} 
 \tilde\alpha^{(2)}\comp\tr\tilde\alpha^{(1)}|\chowhom_m(Y)&=0 \label{eqn:obs1}\\
\tilde\alpha^{(1)}\comp\tr\tilde\alpha^{(2)}|\chowhom_{m-1}(B)&=0. \label{eqn:obs2} 
\end{eqnarray}
\end{observation}
Hypothesis (1)  of Theorem~\ref{result1}    guarantees that \eqref{eqn:obs1} holds.
Indeed, let $B$ be the base locus of the pencil of hyperplane section on $Y$. One has the inclusions
$$
B\into   Y=W\cap\bP^{N+1}\into W
$$
Consider the $(k+1)$-plane $P$ of codimension $k+1$ linear sections sections of
$W$ through $B$. This gives a family of $n$-folds
parametrized by $P$, all of which are in the same family as
$X$:
\[
\begin{array}{ccc}
\alpha_P/P &\longrightarrow & \blow_B W\\
\mapdown{}&& \\
C_P/P&  & \\
\end{array}
\]
and hence
\[
\begin{array}{ccc}
\chowhom_{m}(W) &\mapright{\tilde\alpha_P^{(2)}
\comp\tr\tilde\alpha_P^{(1)}} &\chowhom_{m-1}(B)\\
\mapup{}&&\Vline{}{}{2.7ex}\\
\chowhom_m(Y) &\mapright{\tilde\alpha^{(2)}
\comp\tr\tilde\alpha^{(1)}}&\chowhom_{m-1}(B).
\end{array}
\]
 
 This  diagram together with $\chowhom_m(W)=0$ implies that  $ \tilde\alpha^{(2)}\comp\tr\tilde\alpha^{(1)}=0$.

Now \eqref{eqn:obs2} follows from \eqref{eqn:obs1}, using the other assumptions. Indeed, since  $\chowhom_0(Y)=0$, by Cor.~\ref{diag} there is a
smooth variety $V$ of dimension $2m$, and a $ \beta\in\corr_{-1}(Y,V)$ such that $\tr\beta\comp\beta|\chowhom_m(Y)=\id$.   Consider the correspondence
$\alpha=\tilde\alpha^{(2)}
\comp\tr\tilde\alpha^{(1)}\comp\tr\beta\in\corr_0(V,B)$. This correspondence acts trivially on $\chowhom_{m-1}(V)$ by assumption. Hence, by the reflexive vanishing  property for $B$,  the transpose $\beta\comp\alpha^{(1)}\comp\tr\tilde\alpha^{(2)}$ acts trivially on
$\chowhom_{m-1}(B)$.  One
concludes that $\tilde\alpha^{(1)}\comp\tr\tilde\alpha^{(2)}|\chowhom_{m-1}(B)=0$.
\qed\endproof
  
\section{Complete intersections with small Chow groups}
\subsection{Complete intersections of very small multidegree} \label{subsec:cis}

For complete intersections $X$ one proceeds as follows.  As explained at the start of \S~\ref{sec:extend}, given any linear subspace $L$ contained in $X$,   a generic  complete intersections of sufficiently small multidegree Lemma~\ref{cutout} guarantees the Assumption~\ref{assu}.  By Proposition~\ref{ComplInts} one has $\chow_k(X)=\chowlin_k(X)$, $k\le m$  for sufficiently small multidegrees. The same multidegrees work also for $Y$. So Theorem~\ref{result1} can be applied. Moreover at the same time for very small multidegrees conditions (1) and (2) of this Theorem  will be  verified.   So one is left with conditions (3) and (4).  
 
In the special case $m=1$, i.e. $X$ is a threefold, one needs to see that $\chowhom_2(Y)=0$. This follows from  Cor.~\ref{conseqfor4folds} applied to $Y$ and the fact that $H^3(Y)=0$ by the Lefschetz hyperplane theorem. The latter also implies $H^1(B)=0$  and $b_2(X)=1$. So, for complete intersection threefolds the situation further simplifies:
 
  \begin{cor} \label{for3folds}Let $Y$ be a generic  smooth complete intersection $4$-fold and let $X$   be a smooth hyperplane section. Suppose  that 
\begin{enumerate}
\item $\chowhom_{0}(Y)=0$; 
\item     {\rm GHC}$_1(X)$ hold;
\item $h^{3,0}(X)=0$.
 \end{enumerate}
Then $\chowhom_{0}(X)=0$. 
\end{cor}
It is easy to see that this does not give new cases since Proposition~\ref{ComplInts} for $0$-cycles is Roitman's theorem which is optimal. For weighted complete intersections the situation is more fortunate:
 
 \subsection{Weighted complete intersections} \label{subsec:wci}

For a complete intersection of multidegree $(d_1,\dots,d_r)$ in
$\bP(a_0,\dots,a_{n+r})$ the canonical bundle $K$ is given by the formula
$K=\OO(\sum d_j-\sum a_i)$ and so $X$ is of general type as soon as $\sum
d_j-\sum a_i>0$, even if the canonical line bundle has no sections. For
\emph{quasi-smooth} complete intersections Steenbrink has shown \cite{St} that the cohomology
groups carry pure Hodge structures so that there is an intermediate
Jacobian, an Abel-Jacobi map etc.  Moreover, since such a complete intersection
is the quotient of a smooth complete intersection in $\bP^{n+r}$, working with
rational coefficients, there   is a good intersection theory for the Chow
groups (see \cite{Fu}). Proposition \ref{lowchow} shows that the various Chow groups vanish  for quasi-smooth complete intersections of very small degrees. 

I expect that Cor.~\ref{for3folds} also holds in the weighted situation.  For the moment, I can at least show that everything works as it should in the following 

\begin{example} \rm See \cite[\S~2.9-2.10]{Reid} for the geometry of  this  example. Let $X=X_{6,6,6}$ be a generic complete intersection of type
$(6,6,6)$ in $P=\bP(2^4,3^3)$. One can take for $Y=Y_{6,6}$ a complete intersection of type $(6,6)$ in $P$ and $W=P$. 

The system $|\OO(6)|$ embeds $P$  in  $\bP^{25}$ as the linear join of the Veronese
embedded $\bP^3$, say $v_3(\bP^3)\subset P(\OO(3))=\bP^{19}$ and $\bP^2$, say $v_2(\bP^2) \subset P(\OO(2))=\bP^{5}$. This description makes it quite easy to verify the (appropriate modification of) Assumption~\ref{assu} in this case:  one  just has to makes judicious choices for $Y$, $X'$, $Y'$ as follows.  Any point  $p\in P\subset \bP^{25}$ is on a line $L$  joining two points $q=v_3(a)$ and $r=v_2(b)$. If $p\in X=P\cap H_1\cap H_2\cap H_3$, for generic choices of $H_i$ the point $r$ does not belong to $X$. The hyperplanes containing  $H_1\cap H_2\cap H_3$ form a $\bP^2$ and so one can find two independent hyperplanes $H_1'$ and $H_2'$ in this web passing through $r$. For generic choices $H_1'\cap H_2'$ will meet $P$ transversally and hence $Y=P\cap H_1'\cap H_2'$ is quasi-smooth. By Remark~\ref{cutout2} the general hyperplane $H_3'$ through $L$ will also meet $Y$ transversally and this gives $X'$. Taking a generic  point $s\in v_2(\bP^2)$ for instance, the join with $L$ is a plane $M\subset \bP$ and one can find two independent hyperplanes $H_1'$, $ H_2''$  in the web containing $H_1'\cap H_2'\cap H'_3$ that pass through $s$. These contain $M$ and, again by Remark~\ref{cutout2}  for generic choices   $H_1''\cap H_2''$ meets  $P$ transversally in a quasi-smooth $Y'$.

The threefold $X$ is of general type since  $K_X=\OO_X(1)$ but  $K_X$ has no sections and so $p_g(X)=0$. The condition GHC$_1(X)$   holds due to Ortland's constructions. In fact, any weighted homogeneous polynomial  of degree $6$ can be written $F= C(Y_1,Y_2,Y_3,Y_4)+Q(Z_1,Z_2,Z_3)$ where $C$ is cubic in the first $4$ variables $Y_j$ and $Q$ is quadratic in the remaining ones, $Z_j$. To $X=V(F_1,F_2,F_3)$ then is associated a web $F_\lambda= \sum_{i=1}^3  \lambda_i F_i$ with $\lambda=(\lambda_1,\lambda_2,\lambda_3)\in\bP^2$. Whence a web  of cubics $C_\lambda$ in $\bP^3$ and conics  $Q_\lambda$ in $\bP^2$.
The conics degenerate into a pair of lines $(L,L')$ along a degree $3$ curve $E\subset\bP^2$ which is smooth for generic $X$  and all cubic surfaces contain $27$ lines $L_k$. The  pairs  $(L,L_k)$, $(L',L_k)$ define weighted planes $P'=\bP(2,2,3,3)\subset P$ entirely contained in some $F_\lambda$ and $P'\cap X$ is a curve of genus $2$. The parameter space of those planes is a curve $B$ which is generically $(54:1)$ onto $E$ and whence a correspondence from $B$ to $X$. Ortland has verified that for generic $X$ the Abel-Jacobi map $J(B)\to J^1(X)$ is onto.

Also,   $K_Y=\OO_Y(-4)$, and hence $Y$ is Fano. By Prop.~\ref{lowchow} one sees that $\chowhom_0(Y)=0$.

One could replace $X$ by any toroidal resolution 
$f:\hat X\to X$. One still has $\chowhom_0(\hat X)=0$, since the isolated
singularities are replaced by trees of weighted projective spaces.
The canonical bundle $f^\ast K_X+$ effective divisor
remains big, and $p_g(\hat X)=0$ as well so that $\hat X$ is a smooth
threefold of general type giving an example where the Bloch conjecture
is true.

\end{example}


\begin{thebibliography}{9999}

\bibitem{B} Barlow, R.: Rational equivalence of zero cycles for some
more surfaces with $p_g=0$, Invent. Math \textbf{79} (1985), 303--308.

\bibitem{B-K-L} Bloch, S., A. Kas, D. Lieberman: Zero-cycles
on surfaces with $p_g=0$, Compositio Math \textbf{33} (1976), 135--145.

\bibitem{B-S} Bloch, S., V. Srinivas: Remarks on correspondences and
algebraic cycles, Amer. J. Math {\bf105} (1983), 1235--1253.


\bibitem{E-L} Esnault, H., M. Levine: Surjectivity of cycle maps
``Journ\'ees de g\'eometrie alg\'ebrique d'Orsay'', France, juillet 20-26,
1992, Ast\'erisque \textbf{218} (1993), 203--226.

\bibitem{E-L-V} Esnault, H., M. Levine, E. Viehweg: Chow groups of
projective varieties of very small degree, Duke Math. J. \textbf{87}(1997), 29--58.

\bibitem{Fu} Fulton, W.:{\sl Intersection Theory}, Ergebnisse der Math.,
3.~Folge, \textbf{2}, Springer Verlag, Berlin, Heidelberg, etc, 1984.


\bibitem{Gr1} Griffiths, P.: 
On the periods of certain rational integrals I, resp. II,  Ann. Math. \textbf{90}
(1969), 460-495  resp. 498--541.

\bibitem{Gr2} Griffiths, P.: On the periods of certain rational
integrals III, Publ. Math.  IHES, \textbf{38} (1970), 125--180. 

\bibitem{I-M} Inose, H., M. Mizukami: Rational equivalence of $0$-cycles
on some surfaces of general type with $p_g=0$, Math. Ann
\textbf{244} (1979), 205--217.

\bibitem{Ja} Jannsen, U.: Motivic sheaves and filtrations on Chow 
groups, Proc. Am. Math. Soc. \textbf{55-1} (1991), 245--302.

\bibitem{La} Laterveer, R.:  Algebraic varieties with small Chow groups, J.
Math. Kyoto Univ.  \textbf{38}  (1998), 673--694.

\bibitem{Le} Leoni, M.: Chow groups of weighted hypersurfaces, C. R. Acad. Sci. Paris \textbf{330} (2000), 1--4. 

\bibitem{Ko} Koll{\'a}r, J.:  Rational curves on algebraic varieties, 
Erg. Math. 3. Folge \textbf{32} (1996), Springer, Berlin etc.
 

 
\bibitem{Mum} Mumford, D.: Rational equivalence of zero cycles of surfaces 
 J. Math. Kyoto \textbf{9} (1969), 195--204.

\bibitem{Pa} Paranjape, K.: Hodge theoretic and cycle--theoretic
connectivity, Ann.~Math. \textbf{140} (1994), 641--660. 

\bibitem{PS} Peters, C. and J. Steenbrink: \textit{Mixed Hodge Theory},  Ergebnisse der Math, Wiss. \textbf{52}, Springer Verlag (2008).

\bibitem{Reid} Reid, M.:{A
young person's guide on canonical singularities} {Proc. Am. Math.
Soc.}{46-I} (1985), 345--414.

 
\bibitem{R1} Roitman, A.: Rational equivalence of zero-cycles, Mat.
Sbornik  \textbf{89}(1972), 569--585 = Math. USSR Sbornik \textbf{18} (1972)
571--588.

\bibitem{R2} Roitman, A.: The torsion of the group of 0-cycles modulo 
rational equivalence, Ann. Math. \textbf{111} (1980), 553--570.

\bibitem{Scho} Schoen, C.: On Hodge structures and non-representablity of Chow groups, Compositio Math. \textbf{88} 285--316 (1993) 

\bibitem{St} Steenbrink, J.: Mixed Hodge Structures on the
vanishing cohmology,in {\sl Real and Complex
Singularities}, Sijthoff-Noordhoff, Groningen,(1977), 525--563.

\bibitem{V} Voisin, C.: Sur les z{\'e}ro-cycles de certaines
hypersurfaces munies d'un automorphisme, Ann. Sc. Norm.
Pisa  \textbf{29} (1992), 473--492.

\bibitem{V2} Voisin, C.:  Hodge theory  and complex algebraic geometry II, Cambridge University Press (2003)
\bibitem{V3} Voisin, C.:   Remarks on filtrations on Chow groups and the Bloch conjecture, Annali di matematica \textbf{183}, 421--438 (2004). 

\bibitem{Zu} Zucker, S.: Generalized intermediate jacobians and
the theorem on normal functions, Inv. Mat \textbf{33}
(1976), 185--222.



\end{thebibliography}
\end{document}